\documentclass{amsart}
\usepackage{graphicx,amssymb,amsfonts,amsmath,amsthm,newlfont}
\usepackage{epsfig}
\newtheorem{thm}{Theorem}[section]
\newtheorem{cor}[thm]{Corollary}
\newtheorem{lem}[thm]{Lemma}
\newtheorem{prop}[thm]{Proposition}

\theoremstyle{definition}
\newtheorem{defn}[thm]{Definition}

\newtheorem{example}[thm]{Example}
\theoremstyle{remark}
\newtheorem{rem}[thm]{Remark}
\numberwithin{equation}{section}

\newcommand{\Z}{\mathbb Z}

\newcommand{\N}{\mathbb N}
\newcommand{\Q}{\mathbb Q}

\newcommand{\uP}{\Phi}

\begin{document}

\title[On arithmetic and asymptotic properties of up-down numbers]{}%
\author{F.\ C.\ S.\  Brown, T.\  M.\ A.\ Fink, K.\ Willbrand}%
\address{45 Rue d'Ulm}
\email{brown@clipper.ens.fr}

\maketitle
\begin{center}
\Large{\textbf{On arithmetic and asymptotic properties of up-down
numbers.}}
\end{center}

\begin{abstract}
Let  $\sigma=(\sigma_1,\ldots, \sigma_N)$, where $\sigma_i =\pm
1$,  and let $C(\sigma)$ denote the number of permutations $\pi$
of $1,2,\ldots, N+1,$ whose  up-down signature $\mathrm{sign}(
\pi(i+1)-\pi(i))=\sigma_i$, for $i=1,\ldots,N$.
 We prove that the
set of all up-down numbers $C(\sigma)$ can be expressed by a
single universal polynomial $\Phi$, whose coefficients are
products of numbers from the Taylor series of the hyperbolic
tangent function. We prove that $\Phi$ is a modified exponential,
and deduce some remarkable congruence properties for the set of
all numbers $C(\sigma)$, for fixed $N$. We prove a concise
upper-bound for $C(\sigma)$, which describes the asymptotic
behaviour of the up-down function $C(\sigma)$  in the limit
$C(\sigma) \ll (N+1)!$.

% for
%signatures which are infrequent.% $C(\sigma)\ll (N+1)!$.

\end{abstract}

\newpage
\maketitle
%%%%%%%%%%%%%%%%%%%%%%%%%%%%%%%%%%%%%%%%%%%%%%%%%%%%%%%%%%%%%%%%%%%%%%%%%%%%%%

\section{Introduction}
Let $N\geq 1$, and let $\pi$ be a permutation of $\{1,2,\ldots,
N+1\}$. The up-down signature of $\pi$ is defined to be  the
sequence $\sigma=(\sigma_1,\ldots, \sigma_N) \in \{ 1,-1\}^N$ of
rises and falls of  $\pi$. More precisely, the up-down signature
$\sigma$ is given by the formula:
$$\sigma_i = \mathrm{sign} \big(\pi(i+1)- \pi(i)\big) \qquad \hbox{ for }
\quad 1\leq i\leq N\ .$$ Let $C(\sigma)$ denote the number of
permutations $\pi$ which have  up-down signature $\sigma$. Some
small values of the up-down numbers $C(\sigma)$ are listed in
table 1 below.

%A straightforward way of characterising a landscape, and one which
%is invariant over the choice of distribution of its constituent
%points, is to consider its pattern of increases and decreases,
%without regard to absolute values or the relative values of
%distant points. From this up-down picture many physically
%important properties of the landscape can be deduced, such as the
%number of minima, the size of the largest basin of attraction and
%the distribution of the size of minima.

%The simplest random landscape is a
%sequence of $N+1$ identically and independently
%distributed random numbers. We connect pairs of consecutive
%points with line segments to form a curve.
%If we assume that the probability that two points
%are identical is negligible, we can label these segments
%increasing ($+$) or decreasing ($-$).
%The $N+1$ points can thus be reduced to an up-down signature $\sigma$:
%a string of $+$s and $-$s of length $N$.
%The sequence $0.2,0.6,0.9,0.5,0.3$, for example, has up-down
%signature $++--$.
%Properties such as the number of ascents and the number and
%breadth of minima
%can be directly read off the signature.
%Because the up-down properties of a sequence of points
%depends on their relative values only, the distribution of signatures
%does not depend on the distribution from which the points are drawn.

\begin{small}
\begin{table}[b!]
\begin{center}
\begin{sf}
\begin{tabular}[t]{|c|c||c|c||c|c||c|c||c|c|}
\hline
  \multicolumn{2}{|c||}{N=1}
& \multicolumn{2}{|c||}{N=2} & \multicolumn{2}{|c||}{N=3} &
\multicolumn{2}{|c||}{N=4}
& \multicolumn{2}{|c|}{N=5} \\
\hline \hline $\sigma$ & $C $ & $\sigma$ & $C $ & $\sigma$ & $C $
& $\sigma$ & $C $ &
$\sigma$ & $C $ \\
\hline
    $\,-\,$ &  1  & $\,--\,$ & 1  & $\,---\,$ & 1  & $\,----\,$ &  1 & $\,-----\,$ &  1  \\
    $+$ &  \,1\,  & $-+$ & \,2\,  & $--+$ & \,3\,  & $---+$ &  4 & $----+$ &  5  \\
        &     & $+-$ & 2  & $-+-$ & 5  & $--+-$ &  9 & $---+-$ & 14  \\
        &     & $++$ & 1  & $-++$ & 3  & $--++$ &  6 & $---++$ & 10  \\
        &     &      &    & $+--$ & 3  & $-+--$ &  9 & $--+--$ & 19  \\
        &     &      &    & $+-+$ & 5  & $-+-+$ & 16 & $--+-+$ & 35  \\
        &     &      &    & $++-$ & 3  & $-++-$ & 11 & $--++-$ & \,26\,  \\
        &     &      &    & $+++$ & 1  & $-+++$ &  4 & $--+++$ & 10  \\
        &     &      &    &       &    & $+---$ &  4 & $-+---$ & 14  \\
        &     &      &    &       &    & $+--+$ & 11 & $-+--+$ & 40  \\
        &     &      &    &       &    & $+-+-$ & 16 & $-+-+-$ & 61  \\
        &     &      &    &       &    & $+-++$ &  9 & $-+-++$ & 35  \\
        &     &      &    &       &    & $++--$ &  6 & $-++--$ & 26  \\
        &     &      &    &       &    & $++-+$ &  9 & $-++-+$ & 40  \\
        &     &      &    &       &    & $+++-$ &  4 & $-+++-$ & 19  \\
        &     &      &    &       &    & $++++$ &  1 & $-++++$ &  5  \\
\hline
\end{tabular}
\end{sf}
\end{center}

\vspace{0.2in} \caption{{\small The number   $C(\sigma)$ of
permutations on $N+1$ letters with given up-down signature
$\sigma$ of length $N$. We write $+$ for $+1$ and $-$ for $-1$.
Since $C(\sigma)$ is symmetric on interchanging $+$ and $-$, only
half of the values $C(\sigma)$ for $N=5$ are shown. The maximum
values $A_N$ in each column are asymptotically equal to $2^{N+3}
\pi^{-(N+2)} (N+1)!$ (see \cite{Andre}).
 %Dividing by
%$(N+1)!$ gives the probability $P(\sigma)$ that a random curve has
%the same signature.
}} \label{updown_table}
\end{table}
\end{small}

%The numbers $C(\sigma)$ have many different combinatorial
%interpretations and occur in various different areas of
%mathematics \cite{Bruijn,  MacMahon, Niven,  Szpiro01, Viennot}.
%For example, they are related to  the dimensions of
%representations of the symmetric group and the
%Littlewood-Richardson rule for the multiplication of Schur
%functors \cite{Foulkes}.  The maximum values $1,2,5,16,61,\ldots$,
%occur for the alternating signatures $+-+-\ldots,$ and

The enumeration of permutations with  given up-down signatures is
a long-standing combinatorial problem initiated by Andr\'e
\cite{Andre}, who computed the number of permutations with the
alternating signature of length $N$: $A_N=C(+ - + - \ldots).$ The
numbers $A_N$ are  %the coefficients in the Taylor expansion of
%$\tan z+ \sec z$ and are
called Euler-Bernoulli updown numbers and are given by the Taylor
expansion of $\tan z+ \sec z$.
  These numbers arose in the study of morsifications in singularity theory by
Arnold \cite{Arnold}, who also proved some surprising arithmetic
properties for them. Many variants  of  these numbers have been
studied extensively  by Carlitz and Carlitz-Scoville (see {\it
e.g.}, \cite{Carlitz1,Carlitz2}).
 The numbers $C(\sigma)$ for arbitrary $\sigma$ can
be regarded as a natural generalisation of the numbers $A_N$, but
are altogether less well-understood. They have been studied  in
various combinatorial contexts \cite{Bruijn, MacMahon, Niven,
Szpiro01, Viennot}, and are related, for example, to the
dimensions of irreducible representations of the symmetric group
via the Littlewood-Richardson rule for the multiplication of Schur
functors \cite{Foulkes}.

% representation theory of the symmetric group
%\cite{Foulkes74}, but are altogether less well-understood.

%For example, they are related to

% are given by the Euler and Bernoulli numbers \cite{Andre}.

Now consider $N+1$ independent and identically distributed random
variables $X_1,\ldots, X_{N+1}$, where the $X_i$ are taken from a
continuous distribution ({\it i.e.}, if $i\neq j$, then
$P(X_i=X_j)=0$). Then the quantity
\begin{equation}\label{H}
P(\sigma) = {C(\sigma)  \over (N+1)!}\
\end{equation}
 is the probability that the random curve $X_1,\ldots, X_{N+1}$ has up-down signature $\sigma$.
Thus another motivation for considering the numbers $C(\sigma)$ is
because of the importance of  one-dimensional random energy
landscapes in  statistical physics. These arise in the study of
spin glasses \cite{Derrida,mezard87}, protein folding
\cite{frau97} and drainage networks \cite{frau97}.
%
% in statistical physics, the
%study of spin-glasses \cite{Derrida,frau97,mezard87}
 The numbers $P(\sigma)$ can also be used to define a test for
randomness, which has been applied very effectively to the study
of genetic microarray data in  biology \cite{Us,Us2}.

 It
is also known how to compute the probability that two  random
curves have the same up-down signature \cite{Mallows85}, and how
to compute the expected values of a random permutation with any
given up-down signature \cite{Szpiro01}.

In this paper, we answer questions about the nature of the whole
up-down  sequence (or distribution) for a given length $N$, {\it
i.e.}, the entire set of up-down numbers $C(\sigma)$ (or
$P(\sigma)$), for $\sigma \in \{1,-1\}^N$. This problem is far
from simple because of the highly discontinuous  nature of the
up-down distribution (see fig. 1 below).
\begin{figure}[b!]
  \begin{center}
    \leavevmode
  \mbox{} \hspace{-1 cm} %\vspace{-0.2in}
  \epsfxsize=9.3cm \epsfbox{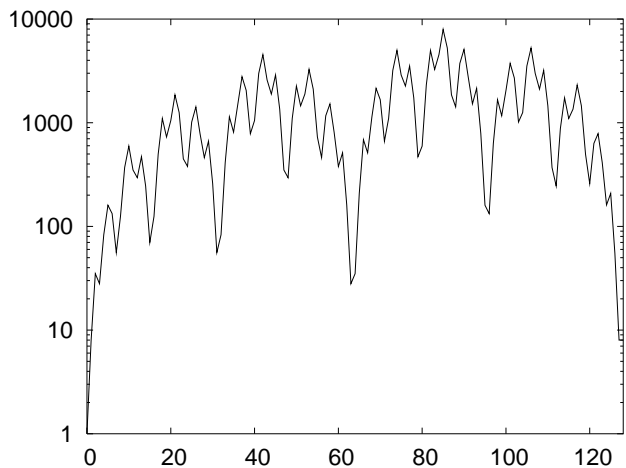}
\caption{
 The number of permutations $C(\sigma)$
as a function of the signature
 $\sigma$ for $N=8$. A number on the
horizontal axis represents a signature via its representation in
binary ({\it e.g.},  for $N=8$, $C(49)=C(00011001)=C(
---++--+)= 1016$). Only the first half of the up-down sequence is
shown. The second half, corresponding to values between 128 and
256, is obtained by symmetry on interchanging $+$ and $-$.
 }
\end{center}
 \label{fractal}
\end{figure}
We approach the problem from two different angles. First of all,
we show that there exists a universal polynomial $\Phi$, whose
coefficients are given by the  Taylor expansion of the hyperbolic
tangent function,  which gives an explicit expression for the
up-down function $C(\sigma)$ for signatures of arbitrary length
(theorem $\ref{thmphi}$). This gives a concise description of the
up-down distribution as the superposition of a small number of
much simpler distributions, and gives an expression for each
up-down number $C(\sigma)$ as an explicit linear combination of
tangent (or Bernoulli) numbers. We also show that the polynomial
$\Phi$ is in fact an exponential with respect to a certain
modified product denoted $\star$ (proposition
$\ref{propphiisexp}$.)
 From
this, we can deduce some remarkable congruence properties
satisfied by the set of numbers $C(\sigma)$ (corollary
$\ref{cor}$). This sheds light on the fine structure of the
distribution $C(\sigma)$.

The second approach is to show  how one can approximate the
up-down distribution $P(\sigma)$  (and hence $C(\sigma)$) by
considering it as a function of the lengths  of its increasing or
decreasing runs. We derive a simple upper-bound for the quantities
$P(\sigma)$ (theorem $\ref{thmapprox}$), which  gives the
asymptotic behaviour of the up-down distribution in the tail
$P(\sigma)\ll 1$. This sheds light on the coarse structure of the
distribution $P(\sigma)$. In applications where the up-down
numbers are used as a test for randomness, this is useful for
establishing the non-randomness of a given data set.

The paper is organised as follows. In $\S2$ we state our main
results.  In $\S3$ we recall some well-known properties of the
up-down numbers, and in $\S4$ we give all  the proofs of  our
results. The second author wishes to thank B. Derrida for some
interesting discussions and suggestions.

\section{Statement of results}

Let
$$\Sigma_N=\{(\sigma_1,\ldots, \sigma_N): \sigma_i \in \{1,
-1\}\}$$
 denote the set of all up-down signatures of length $N$.
Any function $f$ on $\Sigma_N$ can be expressed   as a polynomial
in $N$ variables $s_1,\ldots, s_N$, where $s_i$ takes values in
$\{1,-1 \}$. Since $s_i^2=1$ for all $1\leq i\leq N$, it follows
that $f$ can be written as a sum of linear monomials. For example,
any $\Q$-valued function on the set $\Sigma_2=\{(1,1), (1,-1),
(-1,1), (-1,-1)\}$
%
%
%of four signatures of length two: $++$, $--$, $+-$, $-+$,
 can be
uniquely written in the form:
$$f (s_1,s_2) = a_0 + a_1 s_1 + a_2 s_2 + a_{1,2} s_1 s_2\ ,$$
where $a_0,a_1,a_2,a_{1,2}\in \Q$. Let us define $c_N(s_1,\ldots,
s_N)$ to be the  polynomial function which interpolates the values
of the up-down sequence $C(\sigma)$ for all $\sigma$ of length
$N$. By $(\ref{H})$, the function interpolating $P(\sigma)$ is
given by:
$$p_N(s_1,\ldots, s_N) = {1 \over (N+1)!} c_N(s_1,\ldots, s_N)\ .$$
The first few polynomials $c_1,\ldots, c_5$ are listed below, and
can be used to reproduce all the entries in table $1$.
\begin{eqnarray} \label{qsmall}
c_1 &=& 1\ , \qquad\qquad  c_2 =    \textstyle{1\over 2} (3 - s_1
s_2)\ , \qquad \qquad
c_3 =    3-s_1s_2-s_2s_3\ , \nonumber \\
c_4& =& \textstyle{1\over 2
}(15-5(s_1s_2+s_2s_3+s_3s_4)+2s_1s_2s_3s_4)\ , \nonumber \\
c_5 & = &\textstyle{1\over 2 }(45 -15(s_1s_2+s_2s_3+s_3s_4+s_4s_5)
+ 6(s_1s_2s_3s_4+s_2s_3s_4s_5) + 5s_1s_2s_4s_5) \ .\nonumber
\end{eqnarray}
We will show that the polynomials $c_N$ (and hence $p_N$) can be
obtained by truncating a certain universal polynomial $\Phi$ in an
infinite number of variables $s_1,\ldots, s_N,\ldots.$

\subsection{The universal polynomial.} In order to consider all  up-down sequences
simultaneously,  let
$\Sigma_\infty$ denote the set of  all %infinite sequences
%$(\sigma_1,\ldots, \sigma_n,\ldots)$ such that, for some $N\geq
%1$,  $\sigma_i=0$ for all $i\geq N+1$ and $\sigma_i \in  \{1,-1\}$
%for all $1\leq i\leq N$. In other words
%$$\Sigma_\infty= \bigcup_{N\geq 1} \{(\sigma_1,\ldots,
%\sigma_N,0,\ldots): \,\, (\sigma_1,\ldots, \sigma_N) \in
%\Sigma_N\}$$ is the set of
up-down sequences of arbitrary finite length followed by zeros:
\begin{eqnarray}
\Sigma_\infty= &\{&(\sigma_1,\sigma_2,\ldots, \sigma_n,\ldots):
\hbox{ there exists } N\geq 1 \hbox{ such that } \nonumber \\
&& \sigma_i=0 \hbox{ for all } i\geq N+1\ , \hbox{ and }
\sigma_i\in \{1,-1\} \hbox{ for  all } 1\leq i\leq N \}\ ,
\end{eqnarray}
Let
$$R_N= \Q[s_1,\ldots, s_N] / I_N\ ,$$
where $I_N$ is the ideal generated by the relations
$s_1^2=1,\ldots, s_N^2=1$. Then $R_N$ is naturally identified with
the ring of $\Q$-valued functions on $\Sigma_N$. There are obvious
inclusions $R_N \rightarrow R_{N+1}$ for all $N\geq 1$. The
inductive limit
\begin{equation} \label{induclimit}
R = \lim_{N\rightarrow \infty} R_N \end{equation}
 can naturally be identified with the ring of $\Q$-valued
functions on $\Sigma_\infty$.
 Any element $f\in R$  can be uniquely
written as an infinite series of linear monomials
\begin{equation} \label{monomials}
f(s_1,s_2,\ldots)= a_0+\sum_{k\geq 1} \sum_{0<i_1<\ldots<i_k}
a_{i_1,\ldots, i_k} s_{i_1}\ldots s_{i_k}\ , \quad \hbox{ where }
 a_{i_1,\ldots, i_k}\in \Q\ .\end{equation}
For any $N\geq 1$, we shall write
$$f_N(s_1,\ldots, s_N) = f(s_1,\ldots, s_N,0,0,\ldots)\ ,$$
for the series in $R_N$ obtained by truncating  $f$. The value of
the function $f$ on any signature $\sigma=(\sigma_1,\ldots,
\sigma_N)$ of length $N$ is then given by the finite sum
$f_N(\sigma_1,\ldots, \sigma_N)$.
%can be identified with the ring of infinite series in the
%monomials $s_{i_1}\ldots s_{i_k}$, where $i_1<i_2<\ldots<i_k$.

\begin{defn}
Let $A\subset \N$  denote any  non-empty set of positive whole
numbers, and let  $n=|A|$. The set $A$ can be uniquely partitioned
into $k\geq 1$ maximal runs of $i_1,\ldots, i_k$ consecutive
integers, where $i_1+\ldots +i_k=n$. In other words,
$$A= A_1 \cup \ldots \cup A_k\ , $$
where
$$A_1=\{a_1,a_1+1,\ldots, a_1+i_1-1\} \ , \,\, A_2=\{a_2,a_2+1,\ldots, a_2+i_2-1\} \ ,\,\,\,\ldots \ ,
\,\,
 A_k=\{a_k,a_{k}+1,\ldots a_k+i_k-1\}\ , $$
 such that $a_2 \geq a_1+i_1+1$, \ldots, $a_k\geq
 a_{k-1}+i_{k-1}+1$, so that
 there is a gap between the end of each consecutive
 sequence and the beginning of the next one.
 We call $(i_1,\ldots, i_k)$ the  \emph{run-type} of $A$. The
 run-type of $\{1,2, 4,5 ,7\}$, for example, is $(2,2,1)$.
\end{defn}

\begin{defn} Let $k\geq 1$, and  let $i_1,\ldots, i_k \in \N$.  We define   an infinite series $\gamma(i_1,\ldots, i_k) \in R$
 by the formula:
\begin{equation}
\gamma(i_1,\ldots, i_k)  = \sum_{\emptyset \neq A\subset \N}
\prod_{a\in A} s_a\ ,
\end{equation}
where the sum is over all sets  of positive integers  $A$ which
have  run-type $(i_1,\ldots, i_k)$. The series
$\gamma(i_1,\ldots,i_k)$ is homogeneous of degree $i_1+\ldots
+i_k$.
\end{defn} Writing this out in full gives
\begin{equation}
\gamma(i_1,\ldots,i_k) = \sum (\underbrace{s_{a_1}s_{a_1+1}\ldots
s_{b_1}}_{i_1})(\underbrace{s_{a_2}s_{a_2+1}\ldots s_{b_2}}_{i_2})
\ldots (\underbrace{s_{a_k}s_{a_k+1}\ldots s_{b_k}}_{i_k})\ ,
\end{equation}
where $a_2>b_1+1$, \ldots, $a_k>b_{k-1}+1$.
%where the sum is over all sets of indices $a_1,\ldots,a_k$ subject
%to the conditions
%$$a_2\geq b_1+2\ , \qquad a_3\geq b_2+2 \ ,\qquad \ldots \ , \quad
% a_k\geq b_{k-1}+2\ ,$$
% where $b_1=a_1+2{i_1}-1$, \ldots,
%$b_k=a_k+2{i_k}-1$. In other words, there is a gap between the
%ends of each set of consecutive products which occur in the sum.
%
We have, for example,
\begin{eqnarray}
\gamma(2,2)& =& s_1s_2s_4s_5+
s_1s_2s_5s_6+s_1s_2s_6s_7+s_1s_2s_7s_8+\ldots \nonumber \\
&& + s_2s_3s_5s_6 + s_2s_3s_6s_7+s_2s_3s_7s_8+ \ldots +
s_3s_4s_6s_7 + \ldots+\ldots \nonumber \end{eqnarray} Now consider
the Taylor expansion of the hyperbolic tangent function
\begin{equation}
{\tanh z \over z}= 1 + \sum_{k\geq 1} T_k \,z^{k}\ ,
\end{equation}
where $T_{2k+1}=0$ for all $k$, and $T_2=-1/3$, $T_4=2/15$,
$T_6=-17/315$, $T_8=62/2835$, and, in general,
\begin{equation}
T_{n-2}=  {2^{n} (2^{n}\!-\!1)B_{n}\over n!}\ ,\qquad \hbox{for }
n\geq 4\ ,
\end{equation} where $B_{n}$
is the $n^{\mathrm{th}}$ Bernoulli number.
\begin{defn}
 We define  the \emph{universal polynomial} $\uP\in R$ to be the
 series
\begin{equation}
\label{univ}  \uP=1+ \sum_{\emptyset \neq A \subset \N}
 T_{i_1}\ldots T_{i_k} s_{A_1}\ldots s_{A_k}
\ ,
\end{equation}
where the sum is over all non-empty sets of positive integers $A$,
whose run-type we denote $(i_1,\ldots, i_k)$. As above, the
corresponding partition is denoted  by $A=A_1\cup\ldots \cup A_k$,
and for any non-empty set $B\subset \N$, we write $s_{B} =
\prod_{b\in B} s_b$.
\end{defn}
Equivalently, we can write the universal polynomial in terms of
$\gamma$-series:
\begin{equation}
\label{univ}  \uP=1+ \sum_{k\geq 1} \sum_{i_1\geq 1,\ldots,
i_k\geq 1} T_{i_1}\ldots T_{i_k} \gamma(i_1,\ldots,i_k)\ .
\end{equation}

\begin{thm} \label{thmphi} The universal polynomial describes the up-down
sequences of length $N$ for all $N\geq 1$:
\begin{eqnarray}
p_N(s_1,\ldots,s_N)  & = & 2^{-N} \Phi_N(s_1,\ldots, s_N)\ , \\
c_N(s_1,\ldots,s_N)  & = & (N+1)! \, 2^{-N} \Phi_N(s_1,\ldots,
s_N)\ .
\end{eqnarray}
\end{thm}
\noindent
 Therefore, if
$\sigma=(\sigma_1,\ldots, \sigma_N)$ is any signature of length
$N$, with $\sigma_i \in \{1,-1\}$, then
$$
C(\sigma) = (N+1)!\, 2^{-N} \uP_N(\sigma_1,\ldots, \sigma_N)\ .
$$

%We therefore define two sets of variant polynomials for all $N\geq
%1$:

%By equation $(\ref{H})$, we deduce the following corollary.
%\begin{cor}\label{corforC} For all up-down signatures
%$\sigma=(\sigma_1,\ldots, \sigma_N)$ with $\sigma_i \in
%\{-1,+1\}$,
%\begin{equation}  C(\sigma_1,\ldots, \sigma_N) =
%c_N(\sigma_1,\ldots, \sigma_N)\ . \end{equation}
%\end{cor}
\begin{example}
Consider the case  $N=8$. The corresponding up-down sequence is a
function on the 256 possible up-down signatures of length eight
(see fig. 1 in the introduction). One would expect  the polynomial
in $s_1,\ldots, s_8$ which fits this complicated sequence to have
a large number
 of terms. We have:
  %Thus $q_4 = {1\over 16}(1-
%{1\over 3}\gamma(2) +{2\over 15}\gamma(4))$, ({\it c.f.}
%($\ref{qsmall}$)) and
%$256 q_8$
%consists of the first nine terms of (\ref{tanexpan}).
\begin{eqnarray}\label{phieightexplicit}
\Phi_8 &= &1  - \textstyle{{1\over 3}\gamma_8(2)  +  {2\over
15}\gamma_8(4) + {1\over 9}\gamma_8(2,2)  -
{17\over 315}\gamma_8(6) }  \\
&-&\textstyle{ {2\over 45}(\gamma_8(2,4) +\gamma_8(4,2)) -
\textstyle{{1\over 27}\gamma_8(2,2,2)  + {62\over 2835}
\gamma_8(8)}}\ , \nonumber
\end{eqnarray}
where \begin{eqnarray} \label{gammalistforeight} \gamma_8(2)& =&
s_1s_2 + s_2 s_3 +s_3s_4 + s_4s_5+s_5s_6+s_6s_7 + s_7
s_8 \ , \\
\gamma_8(4) & = & s_1s_2s_3s_4+ s_2s_3s_4s_5 + s_3s_4s_5s_6+
s_4s_5s_6s_7 +s_5 s_6 s_7
s_8 \ ,\nonumber \\
\gamma_8(2,2) & = & s_1s_2 s_4s_5 + s_1 s_2 s_5 s_6 +s_1s_2s_6s_7
+s_1s_2s_7s_8
 + s_2s_3s_5s_6 \nonumber \\
 & & +s_2s_3s_6s_7+s_2s_3s_7s_8 + s_3s_4s_6s_7+s_3s_4s_7s_8+  s_4 s_5 s_7 s_8
 \ ,      \nonumber \\
\gamma_8(6) & = & s_1s_2s_3s_4s_5s_6 + s_2s_3s_4s_5s_6 s_7 + s_3
s_4
s_5 s_6 s_7 s_8\ , \nonumber \\
 \gamma_8(2,4) & = & s_1s_2s_4s_5s_6s_7+ s_1s_2
s_5s_6s_7s_8 + s_2s_3
s_5s_6s_7s_8 \ ,\nonumber \\
 \gamma_8(4,2) & = & s_1s_2s_3s_4s_6s_7+
s_1s_2s_3s_4s_7s_8 +
  s_2s_3
s_4s_5s_7s_8\ ,  \nonumber \\
\gamma_8(2,2,2)& = &s_1s_2s_4s_5s_7s_8\ , \nonumber \\
\gamma_8(8)& = & s_1 s_2 s_3 s_4 s_5 s_6 s_7 s_8 \ .\nonumber
\end{eqnarray}
Recall that the subscript $8$ means that the infinite series
$\gamma$ are truncated up to $s_8$. Quite remarkably, theorem
$\ref{thmphi}$ predicts that
\begin{eqnarray}\label{clevel8}
2\, c_8 (s_1,\ldots, s_8)= 2835 - 945 \,\gamma_8(2) +
378\,\gamma_8(4) +315\, \gamma_8(2,2) -153\, \gamma_8(6)\nonumber \\
-126\,(\gamma_8(2,4)+\gamma_8(4,2)) -105 \,\gamma_8(2,2,2) + 62\,
\gamma_8(8)\ .
\end{eqnarray}
The up-down distribution  for $N=8$  is therefore completely
described by the superposition of just $8$ simpler distributions
$(\ref{gammalistforeight})$, which encode its symmetry in a subtle
and concise way. By truncating further, we retrieve the
polynomials $c_1,\ldots, c_5$ listed earlier.  Note that since
$T_{2k+1}=0$, only $\gamma$'s with even arguments can occur.

%remarkable that $\Phi_8$ contains a total of only 33 monomials.
\end{example}

One can ask in general how many such $\gamma$'s occur in $c_N$. We
thank the referee for pointing out that  this is given
asymptotically by $\alpha^{N-1}$, where $\alpha\approx 1.3247$ is
the real root of $1+x-x^3$ (\cite{Sloane}, sequence A023434).
 This is still exponential, but is considerably slower than $2^N$.
%\begin{rem}
%The polynomial $\gamma_N(2)$ counts the number of local optima in
%a random curve.
%\end{rem}

\subsection{The universal polynomial as an exponential}
The universal polynomial can be succinctly rewritten  as follows.
Let $T$ denote the $\Q$-vector space which is generated by formal
sums of the linear monomials
$$s_{i_1}s_{i_2}\ldots s_{i_k} \ ,\qquad  \hbox{ where }
\qquad i_1<i_2<\ldots<i_k\ .$$ As remarked earlier, $T$ is
isomorphic to the vector space underlying $R$
$(\ref{induclimit})$. We now define a new product $\star: T\otimes
T\rightarrow T$, which is defined on monomials by the formula
$$ (s_{i_1}s_{i_2}\ldots s_{i_k}) \star (s_{j_1}s_{j_2}\ldots
s_{j_\ell}) =
  \begin{cases}
     & s_{i_1}s_{i_2}\ldots s_{i_k} s_{j_1}s_{j_2}\ldots
s_{j_\ell}\quad \text{ if }\quad   j_1>i_k+1\ ,  \\
     & s_{j_1}s_{j_2}\ldots s_{j_\ell} s_{i_1}s_{i_2}\ldots
s_{i_k}\quad \text{ if } \quad  i_1>j_\ell+1\ ,\\
      & 0  \quad  \text{ otherwise} \ ,
  \end{cases}
$$
and extends in the obvious way to all series in $T$. The product
$\star$ makes $T$ into a commutative algebra with unit $1$. We
have, for example, $s_1s_2\star s_2s_3=0$ but $s_1s_2\star s_4s_5=
s_4s_5\star s_1s_2 = s_1s_2s_4s_5$. We define the exponential map
$\exp_\star: T \rightarrow T$ with respect to the product $\star$
by the formula:
$$\exp_\star(a) = 1 + a + {1\over 2!} (a\star a) + {1\over 3!} (a \star
a \star  a) + \ldots \qquad \hbox{ for all } a \in T \ .$$
\begin{prop} \label{propphiisexp} The universal polynomial $\Phi$ is an exponential:
\begin{eqnarray}\label{phiisexp}
\Phi&= & \exp_\star\Big( \sum_{i\geq 1} T_i\,\gamma(i) \Big)
\\
& =&  \exp_\star\big(  T_2\,\gamma(2) \big)\star \exp_\star\big(
T_4 \,\gamma(4) \big) \star \exp_\star\big( T_6 \,\gamma(6)
\big)\star \ldots \ .\nonumber
\end{eqnarray}
%where $\gamma(2k) = \sum_{i\geq 1} s_i s_{i+1}\ldots s_{i+2k-1}$.
\end{prop}
\begin{proof}
Let $\alpha_1,\ldots, \alpha_i,\ldots  \in T$ such that $\alpha_i
\star \alpha_i=0$ for all $i\geq 1$. It is a simple exercise to
show that
$$\exp_\star \Big(\sum_{i\geq 1} \alpha_i \Big)= 1 + \sum_{i} \alpha_i
+ \sum_{i<j} \alpha_i\star \alpha_j + \sum_{i<j<k} \alpha_i\star
\alpha_j \star\alpha_k + \ldots $$ If we apply this argument to
the sum of the infinite series of monomials $T_2\, s_1s_2$,
$T_2\,s_2s_3$,\ldots, $T_4\, s_1s_2s_3s_4$, $T_4\, s_2s_3s_4s_5$,
and so on (recall that $T_{2k+1}=0$), we deduce that
$$\exp_\star\Big( \sum_{i\geq 1} T_i\,\gamma(i) \Big) = 1 +
\sum_{k\geq 1} \sum_{i_1\geq 1,\ldots, i_k\geq 1} T_{i_1}\ldots
T_{i_k} \gamma(i_1,\ldots,i_k)\ ,$$ which proves identity
 $(\ref{phiisexp})$, as required.
\end{proof}

%We can  define the formal logarithm $L$ of $\Phi$ as follows:
%$$L = \log_\star \Phi = \sum_{i\geq 1} b_i\,\gamma(2i) \ .$$
%\begin{cor} The distribution $N^{-1}\,L_N$ tends to a limit as
%$N\rightarrow \infty$.
%\end{cor}
%\begin{proof}It follows from the definition that
%$${1\over N} L_N = \sum^{\lfloor N/2\rfloor }_{i= 1} b_i \, \Big({s_1\ldots s_{2i} +
%\ldots + s_{N-2i+1}\ldots s_N \over N} \Big)\ .$$ Each term in
%parentheses lies in the interval $[0,1]$ for any $N$, so the limit
%exists.
%\end{proof}
%\begin{rem} One can check that $\exp_\star (N^{-1}L_N) \sim \exp (N^{-1}
%L_N)$ as $N\rightarrow \infty$, which implies that the normalised
%logarithm of the up-down distribution $N^{-1} \log \Phi_N$  has a
%limiting distribution as $N\rightarrow \infty$. This was
%conjectured in \cite{Mallows85}.
%As an example, consider the
%infinite alternating sequence $\sigma=(+-+-+-\ldots)$. Then
%$${1\over N} \log \Phi_N \sim 1-b_1+b_2-\ldots = \tan 1 = \pi\ .$$
%\end{rem}

\subsection{Congruence properties for all up-down sequences of fixed length} The universal polynomial can be used to deduce
a number of surprising congruence properties which are satisfied
by the entire up-down sequence of length $N$, for a fixed $N$. We
give two of the most elegant such congruences.
\begin{cor}\label{cor} Let $p$ be any odd prime. For all signatures $\sigma=(\sigma_1,\ldots, \sigma_{N})$
of length $N=p-1,$
\begin{equation} \label{firstcong} C(\sigma) \equiv \sigma_1\ldots \sigma_{p-1} \pmod
p\ .\end{equation} In particular, $C(\sigma)$ only takes the
values $\pm 1 \pmod p$.
 Likewise, for all signatures $\sigma$ of length
$N=p$,
\begin{equation} \label{secondcong} 2 C(\sigma) \equiv (\sigma_1+\sigma_p)
\sigma_1\ldots \sigma_{p} \pmod p\ ,\end{equation} and therefore
$C(\sigma)$ only takes the values $0,\pm 1 \pmod p$.
\end{cor}
The proof of these identities, given in $\S4.2$, will follow from
the formula $(\ref{phiisexp})$ using well-known congruence
properties of Bernoulli numbers due to Kummer and Clausen-Von
Staudt \cite{numbertheory1}.

\begin{example} Many more congruence properties can be derived from $\Phi$ as
follows. For example, in the case $N=8$, we can reduce
$(\ref{clevel8})$ modulo $9$ and $7$ to give the simple relations:
%We deduce that
%$$2 C(\sigma) \equiv -105 \,\gamma(2,2,2) + 62\, \gamma(8) \pmod 9\ ,$$
 %Reducing modulo $9$ and $7$  implies the simple
%relations:
\begin{eqnarray}
2\, c_8  (s_1,\ldots, s_8)& \equiv & -105\, \gamma_8(2,2,2) + 62
\gamma_8(8) \pmod 9\ ,
\nonumber \\
2\, c_8  (s_1,\ldots, s_8)& \equiv & -153\, \gamma_8(6) + 62
\gamma_8(8) \pmod 7 \ ,\nonumber
\end{eqnarray}
This implies that
\begin{eqnarray}
 C(\sigma) &\equiv& (6\,
\sigma_3\sigma_6+4)\sigma_1\ldots \sigma_8
\pmod 9\ , \nonumber \\
 C(\sigma) &\equiv&
(4(\sigma_1\sigma_2+\sigma_1\sigma_8+\sigma_7\sigma_8)
+3)\sigma_1\ldots \sigma_8 \pmod 7\ , \nonumber
\end{eqnarray}
for all signatures $\sigma=(\sigma_1,\ldots, \sigma_8)$ of length
$8$. It follows, for example, that $C(\sigma)$ can only be
congruent to $\pm 1, \pm 2$ modulo $9$ for all $\sigma$ of length
$8$ (see fig. 2).
\end{example}

 \begin{figure}[h!]
 \begin{center}
   \leavevmode
 \mbox{} %\hspace{-1 cm} %\vspace{-0.2in}
 \epsfxsize=6.2cm \epsfbox{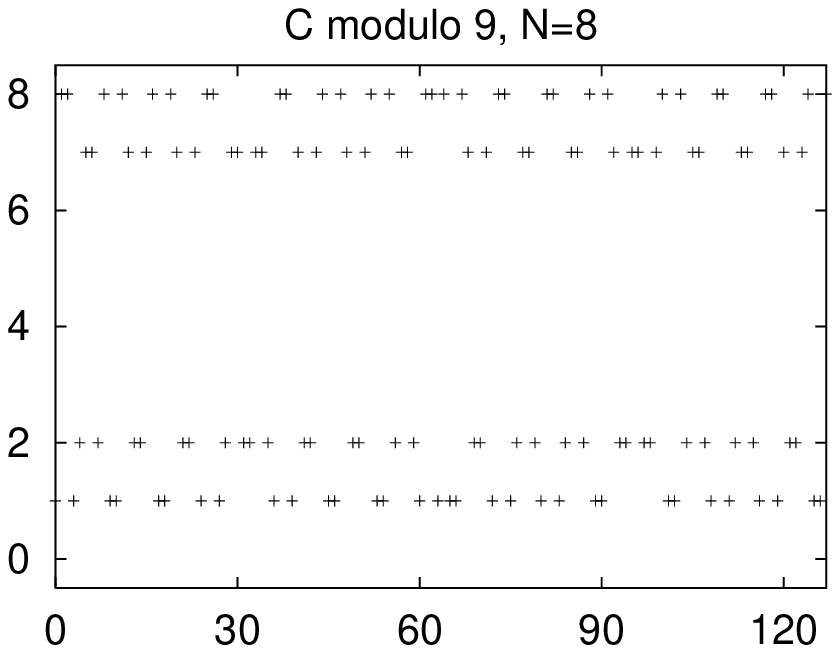}
 \epsfxsize=6.2cm \epsfbox{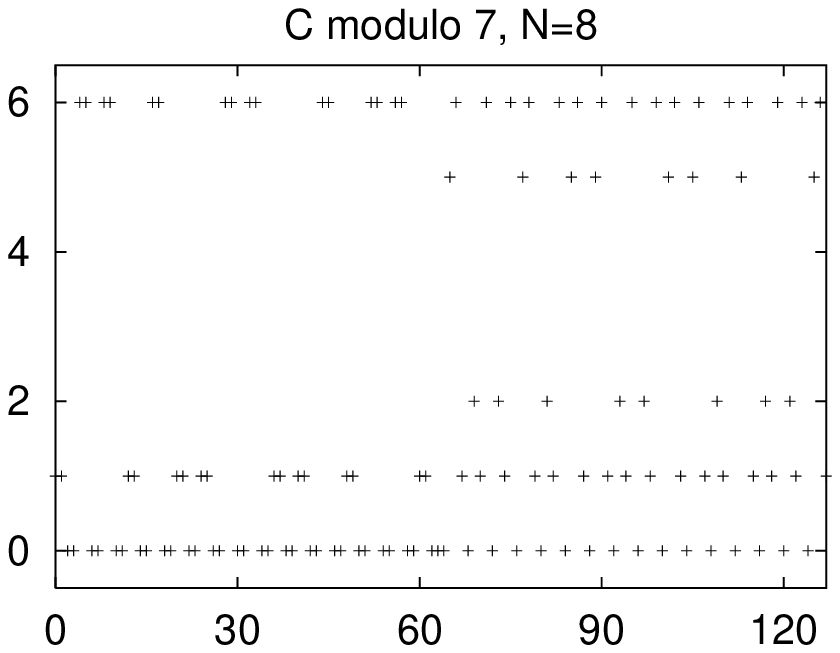}
 \epsfxsize=6.2cm \epsfbox{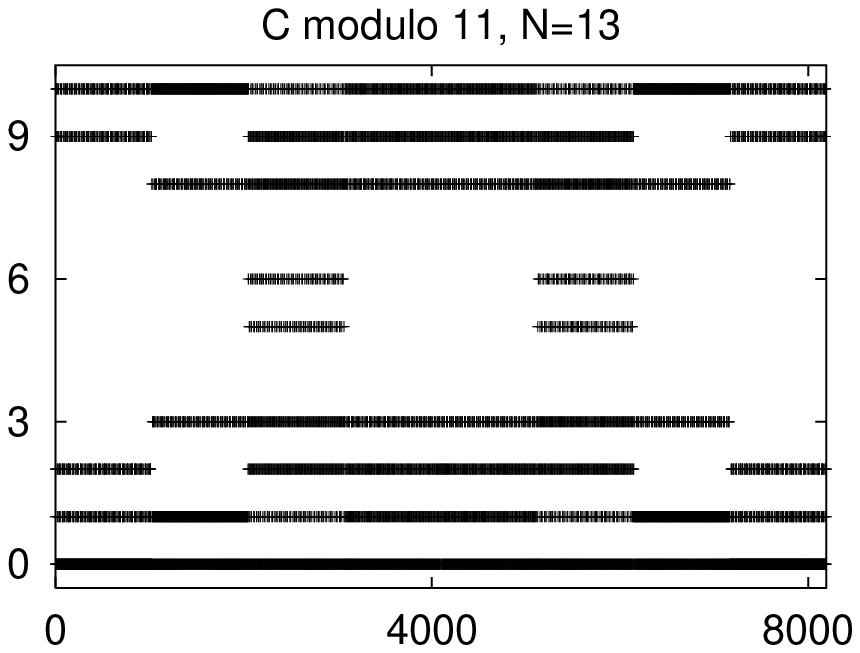}
\epsfxsize=6.2cm \epsfbox{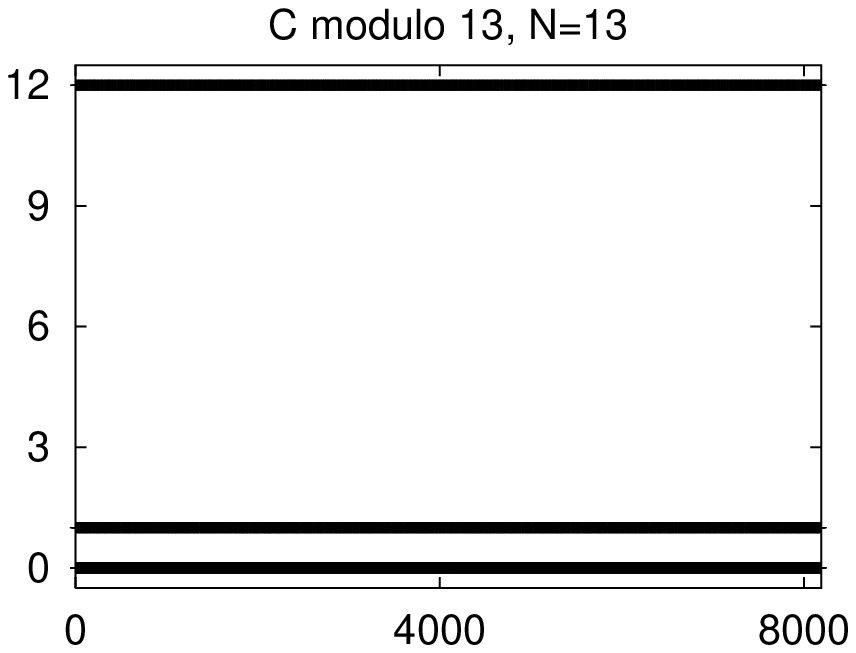}
 \caption{ Congruence properties of the up-down numbers can be
 deduced from the universal polynomial $\Phi$.
  {\sc Top}:
   The (first half of the) up-down numbers $C(\sigma)$ for  signatures $\sigma$ of length $8$, plotted in binary order, modulo 9 (left), and modulo 7
   (right).
% %It follows from the universal polynomial that $C(\sigma) \equiv (6s_3s_6 +4)s \equiv \pm 1, \pm 2 \pmod 9$
% %and   $C(\sigma) \equiv (4(s_1s_2 + s_1s_8 + s_7s_8) +3)s\pmod 7$, where $s=s_1s_2\ldots s_8$.
 %%
 %This shows that the universal polynomial is capturing some  essential features
 %of the up-down distribution.
{\sc Bottom}: A plot of  $C(\sigma)$ modulo $11$ (left) and
     modulo $13$ (right) for all signatures of length $13$.
Equation $(\ref{secondcong})$ predicts that $C(\sigma) \equiv
0,\pm 1 \pmod{ 13}$. The density of the points is such that they
appear as a solid line. } %In general, the up-down numbers for
%signatures of length $N$ satisfy the simplest congruence
%properties modulo primes which are approximately equal to $N$.}
\end{center}
 \label{Congruences}
\end{figure}

\subsection{An upper bound for $P(\sigma)$ and $C(\sigma)$}

% We have seen how
%to compute $P(\sigma)$ for a single signature $\sigma$, or the
%entire set of signatures of a given length.

In  some applications, it is necessary to approximate the
distribution of $C$, or bound $C$ from above. In order to do this,
we need to rewrite a signature $\sigma$ in terms of the lengths of
its runs. Let $(i_1,i_2,\ldots,i_n)$ denote the signature with an
island of $i_1$ pluses, followed by an island of $i_2$ minuses,
and so on, where $i_1+i_2+\ldots+i_n = N$ is a composition of $N$.
For example, $(2,3,1)$ corresponds to $++---+$. One approach to
finding approximations to $C(\sigma)$ is the
 near separability
of  the function $P$ at  maxima or minima. In other words,
probabilistic considerations suggest  the approximation:
%if we divide a
%signature $(i_1,\ldots,i_n)$ into three smaller signatures
%$\rho=(i_1,\ldots,i_{k})$, $\sigma=(i_{k+1},\ldots,i_{\ell})$, $
%\tau=(i_{\ell+1},\ldots, i_n)$, this can be written
%\begin{equation} \label{sep} P(\rho,\sigma,\tau) \simeq P(\rho, \sigma)
%P(\sigma, \tau |\sigma),\end{equation}
\begin{equation} \label{sep} P(i_1,\ldots,i_n)\simeq {P(i_1,\ldots,i_{\ell})
P(i_{k+1},\ldots, i_n) \over P(i_{k+1},\ldots, i_{\ell})}\ ,
\end{equation}
 %where $P(\sigma,\tau
%|\sigma) = P(\sigma,\tau)/P(\sigma)$.
for $k<\ell<n$. Applying $(\ref{sep})$ repeatedly in the case
where $\ell=k+1$, we obtain
%we
%deduce the useful approximation
%\begin{equation} \label{triples}
%P(i_1,.\,.\,,i_n) \simeq { P(i_1,i_2,i_3) P(i_2,i_3,i_4)\ldots
%P(i_{n-2},i_{n-1},i_n) \over P(i_2,i_3)\ldots P(i_{n-2},i_{n-1})}.
%\end{equation}
%This reduces the approximate computation of $P(\sigma)$ for
%arbitrary $\sigma$ to that of $P(i,j)$ (for which there is the
%exact formula $(\ref{ij})$) and $P(i,j,k)$ (whose values could be
%read off a table, for example).
% Many other approximations can be obtained in a similar manner. In
% general, the fewer decompositions  $(\ref{sep})$ used, the better
% the overall approximation.
% If we now apply
%$(\ref{sep})$ repeatedly with $\ell=k+1$, then we obtain a weaker
%approximation, but one which turns out to be an upper bound:
\begin{equation} \label{up1}
P(i_1,.\,.\,,i_n) \simeq { P(i_1,i_2) P(i_2,i_3)\ldots
P(i_{n-1},i_n) \over P(i_2)\ldots P(i_{n-1})}\ .
\end{equation}
 The right-hand side can be written
explicitly in closed form  by $(\ref{ij})$. It turns out that this
approximation is an upper bound, which gives the following
inequality.

\begin{thm} \label{thmapprox} For all $i_1,\ldots, i_n\geq 1$,
\begin{equation} \label{upb}
P(i_1,\ldots ,i_n) \leq {(i_2+1)\ldots(i_{n-1}+1) \over (i_1 + i_2
+ 1) \ldots (i_{n-1} + i_n + 1) } {1\over i_1! \ldots  i_n!}\ .
\end{equation}
\end{thm}
By $(\ref{H})$, we can multiply through by $(i_1+\ldots +i_n+1)!$
 to obtain a similar upper bound  for $C(i_1,\ldots, i_n)$.
\begin{rem} Equation (\ref{upb})
  gives the limiting behaviour of $P$ in the tail
  %end of the   distribution
  $P\ll1$. If the number of islands  $n$ is very small, or if there is a very large island
$i_k$, then certainly the right hand side of $(\ref{upb})$, and
therefore $P(i_1,\ldots,i_n)$ itself, will be small. This  is
relevant when using $P$ as a test for randomness. However, the
converse is far from true,  and the question of when
$P(i_1,\ldots,i_n)$ is small is considerably more subtle. Note
that the denominators $(i_k+i_{k+1}+1)$ in equation $(\ref{upb})$
take into account not just the island sizes $i_k$ but also
first-order dependencies between adjacent islands. One can
speculate that $(\ref{upb})$ is something like the dominant term
of an asymptotic formula expressing $P(i_1,\ldots, i_n)$ in terms
of the island sizes $i_k$.
\end{rem}
\begin{rem}
Equation $(\ref{upb})$ is most accurate when $i_2,\ldots, i_n  \gg
1 $. One can obtain a complementary upper bound  for any
signatures $\rho$ and $\tau$:
\begin{equation}\label{complementary} P(\rho, 1,\tau) \leq
P(\rho)P(\tau).\end{equation} This inequality  follows immediately
from equation $(\ref{L})$. In \cite{Us}, the up-down probabilities
$P(\sigma)$ were used as a test for randomness and applied to
genetic microarray data. By combining the inequalities
$(\ref{complementary})$ and $(\ref{upb})$, one could easily  show
 by hand that many such gene expression curves were non-random.
\end{rem}

\section{Recurrence relations for the up-down numbers}

We  recall two well-known recurrence relations satisfied by the
up-down numbers. The first is linear, the second is quadratic. We
will write, for instance, $(\alpha,j) = (i_1,\ldots,i_n,j)$, where
a Roman letter denotes an island of $+$'s or $-$'s, and a Greek
letter denotes any signature $(i_1,\ldots,i_n)$ which consists of
several islands (see \S2.4).

\subsection{A linear recursion for $C(\sigma)$} The numbers $C(\sigma)$ satisfy the following linear
 recursion relation, which is the same recursion as that for
 multinomial coefficients:
\begin{equation} \label{N4}
C(i_1, \ldots ,i_n) = C(i_1-1,\ldots ,i_n) +
C(i_1,i_2-1,\ldots,i_n) +\ldots+ C(i_1,\ldots ,i_n-1)
\end{equation}
subject to the  boundary conditions $C(0,\alpha)=C(\alpha)$,
$C(\alpha,0) = C(\alpha)$, and
\begin{equation} \label{N5}
C(\alpha,i,0,j,\beta) = C(\alpha, i+j,\beta). %\,\,\, {\rm and}
%\,\,\, C(0,\alpha) = C(\alpha).
\end{equation}
Equation $(\ref{N4})$ can be derived in the  following way (see
also \cite{Carlitz1}). In a permutation of $1,2,\ldots,N+1$ with
signature $(i_1,\ldots,i_n)$, the largest element, $N+1$, must
occur at the end of a sequence of pluses. If we remove it, we
obtain  a permutation of length $N$ with signature $(i_1,\ldots
,i_{2k-1}-1,i_{2k},\ldots,)$ or $(i_1,\ldots,i_{2k-1},
i_{2k}-1,\ldots )$. It follows that there is a one to one
correspondence between the set of all permutations with signature
$(i_1,\ldots,i_n)$ and the union of all permutations with
signatures $(i_1-1,\ldots,i_n)$, \ldots, $(i_1,\ldots,i_n-1)$,
which proves $(\ref{N4})$.

Although there is no simple formula for $C(i_1,\ldots,i_n)$  when
$n\geq 3$, one can show (using the previous recurrence relation,
for example) that
$$
C(i) = 1 \quad \mathrm{and} \quad C(i,j) = \binom{i+j}{i}\ ,
$$
Using the fact that $P(\sigma)= C(\sigma)/ (N+1)!$, for all
signatures $\sigma$ of length $N$, we deduce that
\begin{equation}\label{ij}
P(i) ={1\over (i+1)!} \quad \mathrm{and} \quad P(i,j) = {1\over
(i+j+1)} {1\over i!} {1\over j!}\ .
\end{equation}

%C(i) = 1,
%&\qquad& \textstyle P(i) = {1\over (i+1)!}, \nonumber \\
%\textstyle C(i,j) = \Big({i+j \over i}\Big),
%&\qquad & \textstyle P(i,j) = {1\over (i+j+1)} {1\over i!} {1\over j!}
%The general solution to $(\ref{N4})$ is a
%shifted multinomial coefficient. %lies in the fact that equations
%$(\ref{N5})$ are asymmetric for $n\geq 3$.

\subsection{A quadratic relation for $P(\sigma)$.}
 The second recurrence relation  we will require is most simply written in terms
of $P$.  Let $\sigma$ and $\mu$ be arbitrary signatures.  Then
there is the quadratic relation
\begin{equation} \label{L}
P(\sigma)P(\mu) = P(\sigma + \mu) + P(\sigma - \mu),
\end{equation}
where $\sigma+\mu$  denotes the concatenation of the signatures
$\sigma$, $+$ and $\mu$, and $\sigma-\mu$ is the concatenation of
the signatures $\sigma$, $-$ and $\mu$. In order to obtain
$(\ref{L})$, we interpret $P(\sigma)$ as being the probability
that a random curve has signature $\sigma$. The equation holds
because a random curve $X_1,\ldots, X_{N+1}$ decouples into two
independent sections $X_1,\ldots, X_m$ and $X_{m+1},\ldots,
X_{N+1}$ if one makes no assumption about the
relative values of the points $X_m$ and $X_{m+1}$ where the curves join. %(\cite{Niven}).

\begin{rem} By rewriting $(\ref{L})$ in terms of $\Phi_N=2^{N} p_N$, and
considering the special case when $\mu=\emptyset$, we obtain the
identity
  $$\Phi_N(\sigma) = {1\over 2} \big(\Phi_{N+1}(\sigma+) + \Phi_{N+1}(\sigma-) \big)\ ,$$
 for any signature $\sigma$ of length $N$. This identity
  implies a   self-similarity for the  scaled up-down curves $\Phi_N$:  the values of the up-down sequence of level $N$
   are given by  the average of adjacent values  of the up-down sequence of level $N+1$.
 \end{rem}

\begin{lem} The quantity $P(i_1,\ldots, i_n)$ is given by the  exact formula
$$
%\begin{equation} \label{exact}
%P(i_1,\ldots,i_n)=
\sum_{r_n=0}^{i_n} \sum_{r_{n-1}=0}^{i_{n-1}+r_n} \ldots
\sum_{r_2=0}^{i_2+r_3} {(-1)^{r_2+r_3+\ldots + r_n} \over
(i_n-r_n)!\,(i_{n-1}+r_n-r_{n-1})!\ldots (i_2+r_3-r_2)!
(i_1+r_2+1)!}\ .
$$
%\end{equation}
%
%where the sum is over all integers $r_1,\ldots,r_n\geq 0$ such
%that
%%
%\begin{eqnarray}
%0\,\leq& r_n&\,\leq\, i_n\ , \nonumber \\
%0\,\leq& r_{n-1}+r_n&\,\leq\, i_{n-1}+i_n\ , \nonumber \\
%\vdots\,\, &     & \,\, \vdots \nonumber \\
% 0\,\leq\,&\, r_2+\ldots+r_n\,&\,\leq\, i_2+\ldots +i_n\ , \nonumber
%\end{eqnarray}
%and $r_1$ is defined by $r_1= (i_1+\ldots +i_n +1)-(r_2+\ldots
%+r_n)$.
\end{lem}
\begin{proof} Let $\alpha$ denote any signature, and let $j,k\geq
1$. Applying  equation $(\ref{L})$ with $\sigma=(\alpha,j)$,
$\mu=(k-1)$ implies that $P(\alpha, j, k ) + P(\alpha, j+1,k-1) =
P(\alpha,j)P(k-1)$. Applying this formula inductively and writing
$P(k-1)=1/k!$, we obtain:
\begin{equation}
P(\alpha,j,k) = \sum_{r=0}^k (-1)^r {P(\alpha,j+r) \over (k-r)!}\
.
\end{equation}
This expresses the $P$-value  of an arbitrary signature  in terms
of $P$-values of  signatures  which have a strictly smaller number
of islands. Applying this formula inductively to the signature
$(i_1,\ldots ,i_n)$, one obtains the formula in the lemma.
\end{proof}
%
%
%\begin{eqnarray}
%0\,&\leq\, r_n\,\leq\, i_n, \quad 0&\leq r_n+r_{n-1}\leq i_n+i_{n-1}, \ldots \\
% 0\,&\leq\,\, r_1+\ldots+r_n\,\,&\leq\, i_1+\ldots +i_n. \nonumber
%\end{eqnarray}
%
% This gives an exact formula for $C$ in terms of
%shifted multinomial coefficients, as was expected from $(\ref{N})$.

\begin{rem} Using $(\ref{H})$, the lemma gives an exact formula for
$C(\sigma)$ in terms of multinomial coefficients,
% (as was to be
%expected from the Pascal's pyramid recursion relation
%$(\ref{N4})$).
but which has the disadvantage of being  inefficient to compute. A
similar formula is given in \cite{Szpiro01}, equation $(6)$.
%for computing $P(\sigma)$, however. % (\cite{Carlitz1}). A
%similar formula was obtained in \cite{Szpiro01} by iterated
%integration, and
There are other known methods for computing $C(\sigma)$. For
example, one can express $C(\sigma)$ as the determinant of a
matrix consisting of binomial coefficients (see \cite{Niven},
\cite{MacMahon} and the refinement in \cite{Foulkes}). There is
also
 a simple iterative
algorithm for computing $C(\sigma)$ as a sum of  numbers which are
all positive  \cite{Viennot,Bruijn}, but, unlike the formula given
in the lemma, this does not give a formula in closed form. The
universal polynomial $\Phi$ gives a completely different way to
compute the up-down numbers $C(\sigma)$.
%simply in terms of the computation of the numbers $C(\sigma)$ and
%$P(\sigma)$.
% though the terms in
%the sum $(\ref{exact})$ are all of similar magnitude.
\end{rem}

\section{Proofs}

\subsection{Proof of theorem $\ref{thmphi}$.}
Let $E_{N+1}$ denote the number of permutations on $N+1$ letters
which have an even number of rises. By symmetry, this is also the
number of permutations with an even number of falls. This in turn
is equal to the number of permutations  whose up-down signature
$(\sigma_1,\ldots, \sigma_N)$ satisfies $\sigma_1\ldots
\sigma_N=1$.
\begin{lem} \label{evenrises}We have:
$$E_{N+1} =
   {(N+1)!\over 2}(1+T_N)\ .
$$
\end{lem}
\begin{proof}
Let $A_{n,r}$ denote the number of permutations on $n$ letters
with $r$ rises, where $1\leq r\leq n$. It is well-known
\cite{Carlitz2} that the quantities $A_{n,r}$ are the Eulerian
numbers, whose generating series is given by
$$F(x,y)={e^{xy}-e^x \over ye^x - e^{xy}} =  \sum_{n=1}^\infty {x^n \over n!} \sum_{r=0}^n A_{n,r}
y^r\ .$$ The generating series for permutations with an even
number of rises is therefore given by
$$\sum_{n=1}^\infty E_{n} {x^n \over n!} = {1\over 2} \big(  F(x,1) + F(x,-1) \big) = {1\over 2} \big( {x
\over 1-x} + \tanh(x)\big)\ ,$$ where $F(x,1)$ is to be
interpreted as $\lim_{y\rightarrow 1} F(x,y)$. Comparing the
coefficients of $x^{N+1}/(N+1)!$ yields
%$$i\, \tan ix =  { e^{2x}-1 \over 1 +
%e^{2x}} = {2\over 1+e^{2x}} -1\ .$$
$$E_{N+1}= {(N+1)! \over 2} (1+T_N) \ .$$
\end{proof}

\begin{lem} \label{lemPhimult} Let $N\geq 2$. For all $1\leq n\leq N$,
$$\Phi_N(s_1,\ldots, s_{n-1}, 0, s_{n+1},\ldots, s_N) =
\Phi_{n-1}(s_1,\ldots, s_{n-1}) \Phi_{N-{n}}( s_{n+1},\ldots,
s_N)\ ,$$ where $\Phi_0=1$.
\end{lem}
\begin{proof} If we work in the algebra $T$, the  exponential formula  (proposition $\ref{propphiisexp}$)  gives:
$$\Phi(s_1,\ldots, s_{n-1}, 0, s_{n+1},\ldots ) = \exp_\star \Big(
\sum_{i\geq 1} T_i\, \gamma(i)(s_1,\ldots,
s_{n-1},0,s_{n+1},\ldots ) \Big)\ .$$ By definition of the sums
$\gamma(i)(s_1,s_2,\ldots) = \sum_{k\geq 1} s_k s_{k+1}\ldots
s_{k+i-1}$, this is:
$$ \exp_\star \Big( \sum_{i\geq 1} T_i\, \gamma(i) (s_1,\ldots, s_{n-1},0,0,\ldots)
+ \sum_{i\geq 1} T_i\, \gamma(i) (s_{n+1},s_{n+2},\ldots )\Big)\ .
$$ By the multiplicativity of the exponential, this is a product:
$$ \exp_\star \Big( \sum_{i\geq 1} T_i\, \gamma(i) (s_1,\ldots, s_{n-1},0,0,\ldots)
\Big)\star \exp_\star \Big( \sum_{i\geq 1} T_i\, \gamma(i)
(s_{n+1},s_{n+2},\ldots )\Big)$$
$$= \Phi(s_1,\ldots, s_{n-1},0,0,\ldots) \star
\Phi(s_{n+1},s_{n+2},\ldots)\ .$$ We have proved that
$$\Phi(s_1,\ldots, s_{n-1}, 0, s_{n+1},\ldots )=
\Phi(s_1,\ldots, s_{n-1},0,0,\ldots) \star
\Phi(s_{n+1},s_{n+2},\ldots)\ ,$$ in the algebra $T$.
 But the definition of the product
$\star$ coincides with the ordinary product for monomials which
are sufficiently far apart:
$$s_{i_1}\ldots s_{i_r}\star s_{j_1}\ldots s_{j_k} =
s_{i_1}\ldots s_{i_r}s_{j_1}\ldots s_{j_k}$$ if
$i_1<\ldots<i_r\leq n-1$ and $n+1\leq j_1<\ldots< j_k$. It follows
that the identity
$$\Phi(s_1,\ldots, s_{n-1}, 0, s_{n+1},\ldots ) = \Phi(s_1,\ldots, s_{n-1},0,0,\ldots)
\Phi(s_{n+1},s_{n+2},\ldots)$$ holds in the algebra  $R$. The
result follows on truncating. The lemma can also be proved by
direct computation using the definition of the universal
polynomial  (equation $(\ref{univ})$).
\end{proof}

{\em Proof of theorem $\ref{thmphi}$.} For all  $N\geq 1$, there
exists a polynomial $p_N(s_1,\ldots, s_N)\in R_N$ such that
$P(\sigma) = p_N(\sigma_1,\ldots, \sigma_N)$ for all signatures
$\sigma=(\sigma_1,\ldots, \sigma_N)$.  We can write $p_N$ uniquely
as  a linear monomial in $s_1,\ldots, s_N$ with coefficients in
$\Q$. First of all, the quadratic relation $(\ref{L})$ implies
that
\begin{eqnarray} p_{n-1} (s_1,\ldots,
s_{n-1})\,p_{N-n}(s_{n+1},\ldots,s_N) & = & p_N(s_1,\ldots,
s_{n-1},1, s_{n+1},\ldots, s_N) \nonumber \\ & & + \quad
p_N(s_1,\ldots, s_{n-1},-1, s_{n+1},\ldots, s_N)\ , \nonumber
\end{eqnarray} for all $1\leq n\leq N$. This can  be rewritten:
\begin{equation}\label{qprodid}
2\, p_N(s_1,\ldots, s_{n-1},0, s_{n+1},\ldots, s_N) =
 p_{n-1}
(s_1,\ldots, s_{n-1})\,p_{N-n}(s_{n+1},\ldots,s_N)
\end{equation}
 Suppose by induction that $p_n=2^{-n}\, \Phi_n$ for all $1\leq
n< N$. Then lemma $\ref{lemPhimult}$ implies that the polynomial
$2^{-N} \Phi_N$ satisfies  identity $(\ref{qprodid})$ also.
 It
follows from the  induction hypothesis that $p_N$ and $2^{-N}\,
\Phi_N$ coincide whenever at least one of the $s_i$'s is $0$.
Since only linear monomials are involved, this implies that
$p_N-2^{-N} \Phi_N$ is a multiple of $s_1\ldots s_N$. In order to
compute the coefficient of the term $s_1\ldots s_N$, let
$$S= \{(s_1,\ldots,s_N) :\,\, s_i\in \{\pm 1\} \hbox{ such that }
s_1\ldots s_N=1\}\ .$$ For any $ 1\leq i_1<\ldots <i_k\leq N$,
where  $k$ is strictly smaller than $N$, we have:
$$\sum_{(s_1,\ldots, s_N) \in S} s_{i_1}\ldots s_{i_k} = 0 \ .$$
It follows that taking the sum over all signatures in $S$ picks
out the constant term $1$ and the leading term $s_1\ldots s_N$
only. It therefore suffices to show that
\begin{equation}\label{finalstep}
\sum_{(s_1,\ldots,s_N)\in S} p_N(s_1,\ldots, s_N)=
\sum_{(s_1,\ldots,s_N)\in S}2^{-N} \Phi_N(s_1,\ldots, s_N)\ .
\end{equation}
 The left
hand side is the probability that the signature
$\sigma=(\sigma_1,\ldots, \sigma_N)$ of a random curve satisfies
$\sigma_1\ldots \sigma_N=1$. This is just $E_{N+1}/(N+1)!$, where
$E_{N+1}$ is the number of permutations on $N+1$ letters which
have an even number of rises. The right hand side is
$$ |S|\, 2^{-N} (1 + \hbox{coeff. of } s_1\ldots s_N \hbox{ in } \Phi_N) =
   2^{-1}(1+T_N) \ .
$$
%We can write the polynomial $p_N$ as a sum of monomials $p_0
%+\ldots + p^{max} s_1\ldots s_N$,  where $p_0=2^{-N}$, which
%implies that
%$$|S|(2^{-N}+ p^{max}) = \sum_{ (s_1,\ldots, s_N) \in S}
%p_N(s_1,\ldots, s_N) \ .$$ This quantity is just given by the
%number $E_N$ of permutations of length $N+1$ with an even number
%of rises. Similarl
By lemma $\ref{evenrises}$, both sides of $(\ref{finalstep})$
agree, which completes the induction step. We conclude that $p_N=
2^{-N} \Phi_N$, as required.
 \hfill{$\square$}

\subsection{Proof of corollary $\ref{cor}$.}
First recall the theorem due to
 Clausen-Von Staudt  (\cite{numbertheory1}, theorem 5.10),
which states that for all $k\geq 2$,
\begin{equation} \label{VonStaudt}
B_{2k} - \sum_{(p-1)|2k}{1\over p}\in \Z \ ,\end{equation} where
the sum ranges over all primes $p$ such that $p-1$ divides $2k$.
Now, the coefficients which occur in  the polynomial
$c_N(s_1,\ldots, s_N)$ are
\begin{equation}\label{cNcoeff}
{(N+1)! \over 2^N} T_{2k-2} = \Big({(N+1)! \over 2^N}\Big)\Big(
{2^{2k} (2^{2k}-1)B_{2k} \over (2k)!}\Big)  \qquad \hbox{ for } 4
\leq 2k \leq N+2\ . \end{equation}
%Bernoulli number $B_{2k}$ are precisely those primes $q$ such that
%$q-1$ divides $2k$.
Now let $p$ be an odd prime, and suppose that $N= p-1$.  If $4\leq
2k\leq N-2$, the prime $p$ does not occur in the denominator of
$B_{2k}$ by $(\ref{VonStaudt})$, and therefore
$${(N+1)! \over 2^N} T_{2k-2} \equiv 0 \pmod p\ .$$
It remains to compute the coefficients $(\ref{cNcoeff})$ for
$2k=N$ and $2k=N+2$. In the first case, we have
$${(N+1)! \over 2^N} T_{N-2} = \Big({(N+1)! \over 2^N}\Big)\Big( {2^{N} (2^{N}-1)B_{N} \over
N!}\Big)=p \,B_{p-1} (2^{p-1}-1)\ .$$ But $2^{p-1}-1 \equiv 0\pmod
p$,  and  $p\,B_{p-1}\equiv 1 \pmod p$ by $(\ref{VonStaudt})$. It
follows that this coefficient vanishes modulo $p$ also. Therefore
all terms in the polynomial $c_N$ vanish modulo $p$ except the
leading term, and we are left with:
$$c_N(s_1,\ldots, s_N) \equiv (N+1)!\, 2^{-N}\,T_{N} \,\gamma_N(N) \pmod p\ ,$$
where $\gamma_N(N)$ consists of the single term $s_1\ldots s_N$.
By equation  $(\ref{cNcoeff})$, we have
 $$ {(N+1)! \over 2^N}
T_{N} = \Big({p! \over 2^{p-1}}\Big)\Big( {2^{p+1}
(2^{p+1}-1)B_{p+1} \over (p+1)!}\Big)\equiv 12 \,B_{p+1} \pmod p\
.
$$
The  congruences for Bernoulli numbers  discovered by Kummer
(\cite{numbertheory1}, corollary 5.14) imply, in particular, that
$2\, B_{p+1}\equiv (p+1) B_2 \pmod p$, and so $12\, B_{p+1} \equiv
1 \pmod p$. We conclude that
$$c_N(s_1,\ldots, s_N) \equiv  \gamma_N(N)=s_1\ldots s_N \pmod p\ ,$$
as required.

The result when  $N= p$ holds for similar reasons, since all the
terms of $c_N$ vanish modulo $p$ except the leading term. The
coefficient of this term is:
 $$ {(N+1)! \over 2^N}
T_{N-1} = \Big({(p+1)! \over 2^{p}}\Big)\Big( {2^{p+1}
(2^{p+1}-1)B_{p+1} \over (p+1)!}\Big)\equiv 6\,  B_{p+1} \equiv
2^{-1} \pmod p\ .
$$
This proves that
$$2\, c_N(s_1,\ldots, s_N) \equiv \gamma_N(N-1) =
s_1\ldots s_{N-1} + s_2\ldots s_{N} \pmod p\ ,$$ as required, and
completes the proof of corollary $\ref{cor}$.

% All terms of $c_N$ vanish modulo $p$
%except the last, by the same argument as above.
%$$= {p+2 \over 2^N} 2^{p+1} (2^{p+1} - 1 ) B_{p+1}\ .$$
%This is congruent to
%$$12\,  B_{p+1}  \equiv (p+1) \,6\, B_2 \equiv 1 \pmod p\ .$$

\subsection{Proof of theorem $\ref{thmapprox}$}

We first prove some general inequalities relating up-down numbers
for different signatures of equal length. A similar-looking
inequality was proved by Niven \cite{Niven} to prove that the
value of $C(\sigma)$ is greatest on the  alternating signature
$\sigma=+-+-\ldots$.
\begin{lem}\label{leminitialcase} Let $\alpha$ denote any signature, and let $a,b,c \in \N$ such that $a\geq c$. Then
$$C(\alpha,a-c+1,b,c) \geq C(\alpha, a+1,b)\ .$$
\end{lem}

\begin{proof}
This inequality is easily proved by induction with respect to the
total length $\ell=|\alpha|+a+b+1$, where $|\alpha|$ is the length
of the signature $\alpha$. The details are left to the reader. The
induction step is given by rewriting the left-hand side using
relation $(\ref{N4})$:
$$C(\alpha,a-c, b, c) + C(\alpha,a-c+1,b-1,c) + C(\alpha,a-c+1,b,c-1)$$
plus terms of the form $C(\alpha', a-c+1,b,c)$, where $\alpha'$ is
a signature of shorter length than $\alpha$. Likewise, the
right-hand side can be written
$$C(\alpha,a,b) + C(\alpha,a+1,b-1)\ ,$$
plus terms of the form $C(\alpha', a+1,b)$. If we assume that the
inequality  holds for  $\alpha$ with all smaller values of
$a,b,c$, then $C(\alpha,a-c+1,b-1,c)\geq C(\alpha,a+1, b-1)$ (this
is the case $(\alpha,a,b-1,c)$),
 and $C(\alpha,a-c+1,b,c-1) \geq C(\alpha,a,b)$ (this is the case
 $(\alpha,a-1,b,c-1)$).
  If we assume that the inequality holds for all
 signatures $\alpha'$ of shorter length than $\alpha$, and
  $a,b,c$, then $C(\alpha', a-c+1,b,c)\geq
 C(\alpha',a+1,b)$. This is enough to complete the induction step,
 and hence the proof. The initial cases  $b=0$, $c=0$ are both trivial by
 $(\ref{N5})$. The case $a=c$ is proved using an inductive argument similar
 to the one given above.
\end{proof}

\begin{prop}
Let $\alpha$ denote any signature, and let $a,b,c\in \N$ such that
$a\geq c\geq 1$. Then for all $0\leq n\leq c-1$,
\begin{equation}\label{claim}
 C(\alpha, a-n,b,c) \geq C(\alpha, a+1, b ,c-n-1)\ .
 \end{equation}

\end{prop}
\begin{proof}
The proof is by induction on the total length
$$\ell(\alpha,a,b,c,n) = |\alpha| + a + b+c -n\ ,$$
where $|\alpha|$ denotes the length of the signature $\alpha$. Let
$a',b',c',n'\in \N$ such that  $a'\geq c'\geq 1$, and $0\leq
 n'\leq c'-1$.  Suppose that $(\ref{claim})$  is true for all :
$$a \leq a'\ , \quad  b\leq b' \ , \quad   c\leq c'\ ,\quad    n \geq
n'\ ,$$ and all $\alpha$ satisfying  $|\alpha|\leq |\alpha'|$ such
that
$$a \geq c \geq 1\quad\ , \quad c-1\geq  n \quad \ ,\quad \hbox{ and }   \quad
\ell(\alpha,a,b,c,n) < \ell(\alpha',a',b',c',n') \ .$$ Then we
will prove $(\ref{claim})$ for $a',b',c'$ and $n'$.
% ({\it i.e.},
%the induction goes  downwards in $n'$, and upwards in $a',b',c'$,
%$|\alpha|$).
 First of all, let us assume that $b'>0$ and
$c'-n'-1>0$. This implies that $a'-n'\geq 1$. By $(\ref{N4})$,
$$C(\alpha,a'-n',b',c')= C(\alpha,a'-n'-1,b',c') + C(\alpha,a'-n',b'-1,c') + C(\alpha,a'-n',b',c'-1)$$
plus terms of the form $C(\alpha',a'-n',b',c')$, where $\alpha'$
is strictly shorter than $\alpha$. Each term in the right-hand
side can be bounded below by the induction hypothesis. The middle
term is bounded below as follows:
\begin{equation}\label{proof1}
C(\alpha,a'-n',b'-1,c') \geq C(\alpha,a'+1,b'-1,c'-n'-1)\ .
\end{equation} Similarly,   on setting  $a=a'-1$, $c=c'-1$, $b=b'$,
$n=n'-1$, we obtain:
$$C(\alpha,(a'-1)-(n'-1),b',c'-1) \geq C(\alpha,(a'-1)+1,b',(c'-1)-(n'-1)-1)$$
{\it i.e.},
\begin{equation}\label{proof2}
C(\alpha,a'-n',b',c'-1) \geq C(\alpha,a',b',c'-n'-1)\ .
\end{equation}
 Finally, we have
$$C(\alpha,a'-(n'+1),b',c') \geq C(\alpha,a'+1,b',c'-(n'+1)-1)\ ,$$
which is just:
\begin{equation}\label{proof3}
C(\alpha,a'-n'-1,b',c') \geq C(\alpha,a'+1,b',c'-n'-2)\ .
\end{equation} Adding the three inequalities $(\ref{proof1})$,
$(\ref{proof2})$ and  $(\ref{proof3})$ together,  we obtain
$$  C(\alpha,a'-n'-1,b',c') + C(\alpha,a'-n',b'-1,c') + C(\alpha,a'-n',b',c'-1)            $$
$$\geq C(\alpha,a',b',c'-n'-1) +C(\alpha,a'+1,b'-1,c'-n'-1) +
 +C(\alpha,a'+1,b',c'-n'-2)$$
 After adding  inequalities of the form
 $C(\alpha', a'-n',b',c') \geq C(\alpha', a'+1,b',c'-n'-1)$,
and rewriting the left and right-hand sides using $(\ref{N4})$, we
obtain
$$C(\alpha,a'-n',b',c') \geq  C(\alpha,a'+1,b',c'-n'-1)$$
which proves $(\ref{claim})$ for $a',b', c' $ and $n'$.

We need to check  the initial cases when  $b=0$, $c=n+1$, or
$|\alpha|=0$. If $b=0$, then $(\ref{claim})$ is trivial, since, by
$(\ref{N5})$, $C(\alpha,a-n,0,c)= C(\alpha, a+c-n)=C(\alpha,
a+1,0,c-n-1)$. If $n=c-1$, then $(\ref{claim})$ reduces to the
inequality of lemma $\ref{leminitialcase}$. The case when
$|\alpha|=0$ clearly holds  from the induction argument given
above. Likewise, the case $a=c$ is also covered by the argument
above.
%If $c=1$ then $n=0$ and we are in the
%situation of the previous lemma also.
%This proves the proposition.
% Finally, observe that if $\alpha=(i_1,\ldots,i_n)$ and $i_n$ hits zero, since,
%after concatenation, this just makes the constant $a$ much bigger,
%so the statement is still true. By the way in which we have done
%the induction, $a$ can never hit zero, since $a\geq c\geq 1$ at
%each step.
\end{proof}

\begin{cor}\label{finalcor} For any signature $\alpha$, and $a,b,c\in \N$ such
that  $a\geq c$, we have
 $$C(\alpha, a,b,c) \geq C(\alpha, a+1, b ,c-1)\ .$$
 Equivalently, $P(\alpha, a,b,c) \geq P(\alpha, a+1, b ,c-1).$
\end{cor}
\begin{rem}The  corollary implies  that
$C(\alpha,i,j,k)$ is maximised (for values of $i\geq k$ such that
$i+k$ is fixed) when $i$ and $k$ are most nearly equal.
\end{rem}
%One should compare niven's inequality. I think it is totally
%unrelated, but check anyway.

\emph{Proof of theorem $\ref{thmapprox}$}. Let $\gamma$ denote any
up-down signature. We write $\gamma=(\beta,r+1)$, where $r\geq 0$.
It is clear that
$$P(\beta,r) P(j+1,k,j) = P(\beta,r) P(j,k,j+1)\ .$$
Using relation $(\ref{L})$, this implies that:
$$P(\beta,r+1,j+1,k,j) + P(\beta,r, j+2,k,j) =
P(\beta,r+1,j,k,j+1) + P(\beta,r, j+1,k,j+1)\ .$$ Corollary
$\ref{finalcor}$ implies that $P(\beta,r, j+2,k,j)\leq P(\beta,r,
j+1,k,j+1) ,$ on setting $\alpha=(\beta,r)$, $a=j+1$, $b=k$, and
$c=j+1$. Substituting into the previous equality implies that
$$ P(\beta,r+1,j+1,k,j)\geq P(\beta,r+1,j,k,j+1)  \ .$$
Recalling that $\gamma= (\beta,r+1)$,   this is just
\begin{equation}
P(\gamma, j+1,k,j)  \geq P(\gamma, j,k,j+1) \ ,\end{equation}
which, by adding $P(\gamma,j,k+1,j)$ to both sides, implies that
$$P(\gamma, j, k+1, j) + P(\gamma, j, k, j+1)\leq
P(\gamma, j+1,k,j)+ P(\gamma, j, k+1, j)\ . $$ By $(\ref{L})$,
this is equivalent to the inequality:
$$ P(\gamma,j,k) P(j) \leq P(\gamma, j) P(k,j) \ . $$
It follows from $(\ref{ij})$ that
$$P(\gamma, j,k) \leq {P(\gamma,j ) P(j,k) \over P(j)} = P(\gamma,j) {j+1 \over j+k+1} {1\over k!}\ .$$
Applying this inequality inductively to the up-down sequence
$(i_1,\ldots,i_n)$, we obtain:
$$P(i_1,\ldots, i_n) \leq  P(i_1,\ldots, i_{n-1}) {i_{n-1}+1 \over i_{n-1}+i_n+1} {1\over   i_n!}\leq \ldots $$
$$ \ldots \leq {(i_2+1)\ldots (i_{n-1}+1) \over
(i_1+i_2+1)\ldots (i_{n-1}+i_n+1)} {1\over i_1!\ldots i_n!} \ ,$$
 which is  precisely inequality $(\ref{upb})$.
\hfill{$\square$}

%----------------------------------------------------------------------------

%\begin{small}

%\end{small}

\end{document}